\documentclass{elsart3-1}

\usepackage[centertags]{amsmath}
\usepackage{amsfonts,amssymb}

\usepackage[english,francais]{babel}

\newtheorem{theorem}{Theorem}[section]

\newtheorem{e-proposition}[theorem]{Proposition}

\newtheorem{e-definition}[theorem]{Definition\rm}

\newtheorem{theoreme}{Th\'eor\`eme}[section]

\newtheorem{proposition}[theoreme]{Proposition}

\newcommand{\LL}{\mathbb{L}}

\renewcommand{\d}{\partial}
\newcommand{\Zf}        {{\mathcal{Z}_*}}
\newcommand{\PP}{\mathbb{P}}
\newcommand{\ZZ}{\mathbb{Z}}
\newcommand{\solf}        {\alpha}
\newcommand{\TT}{\mathbb{T}}
\newcommand{\RR}{\mathbb{R}}

\newcommand{\expect}{\mathbb{E}}
\newcommand{\EE}{\expect}

\newcommand{\sol}        {w}
\newcommand{\base}        {e}

\def\og{\leavevmode\raise.3ex\hbox{$\scriptscriptstyle\langle\!\langle$~}}
\def\fg{\leavevmode\raise.3ex\hbox{~$\!\scriptscriptstyle\,\rangle\!\rangle$}}

\begin{document}

\begin{frontmatter}

\selectlanguage{english}

\title{Malliavin calculus and ergodic properties of highly degenerate 2D stochastic Navier--Stokes equation}

\selectlanguage{english}
\author[authorlabel1]{Martin Hairer}
\ead{hairer@maths.warwick.co.uk}
\author[authorlabel2]{Jonathan C.\ Mattingly}
\ead{jonm@math.duke.edu}
\author[authorlabel3]{\'Etienne Pardoux}
\ead{pardoux@cmi.univ-mrs.fr}

\address[authorlabel1]{University of Warwick, Coventry CV4 7AL, UK}
\address[authorlabel2]{Math Department, Duke University, Box
  90320, Durham, NC 27708  USA}
\address[authorlabel3]{LATP/CMI, Universit\'e de Provence, 39 rue F. Joliot Curie, 13 453 Marseille cedex 13, France}

\begin{abstract}
 The objective of this note is to present the results from the two
  papers \cite{b:MattinglyPardoux03Pre} and
  \cite{b:HairerMattingly04}.  We study the Navier--Stokes equation on
  the two--dimensional torus when forced by a finite dimensional white
  Gaussian noise. We give conditions under which both the law of the
  solution at any time $t>0$, projected on a finite dimensional
  subspace, has a smooth density with respect to Lebesgue measure and
  the solution itself is ergodic.  In particular, our results hold for
  specific choices of four dimensional white Gaussian noise. Under
  additional assumptions, we show that the preceding density is
  everywhere strictly positive.
\vskip 0.5\baselineskip

\selectlanguage{francais}
\noindent{\bf R\'esum\'e}
\vskip 0.5\baselineskip
\noindent
Le but de cette Note est d'annoncer les r\'esultats contenus dans
les articles \cite{b:MattinglyPardoux03Pre} et \cite{b:HairerMattingly04}.
Nous \'etudions l'\'equation de Navier--Stokes sur le tore bidimensionnel 
forc\'ee par un bruit blanc gaussien de
dimension finie. Nous donnons des conditions sous lesquelles 
d'une part la loi de la solution \`a tout instant $t>0$, projet\'ee sur un 
espace de dimension finie, a une densit\'e r\'eguli\`ere
par rapport \`a la mesure de Lebesgue (qui, sous des hypoth\`eses 
suppl\'ementaires, est strictement positive partout), et d'autre part la 
solution de la m\^eme \'equation est un processus ergodique. En particulier 
ces r\'esultats sont vrais dans certains cas de bruit blanc gaussien de 
dimension quatre. 
\end{abstract}
\end{frontmatter}

\section{Introduction}

This note reports on recent progress made in \cite{b:MattinglyPardoux03Pre,b:HairerMattingly04} 
on the study of the two dimensional
Navier--Stokes equation driven by an additive stochastic forcing. Recall
that the Navier--Stokes equation describes the time evolution of an
incompressible fluid. In vorticity form, it is given by 
\begin{equation}
  \label{eq:vorticity}
 \left\{ \begin{aligned}
  &\frac{\partial \sol}{\partial t}(t,x)+  B(\sol,\sol)(t,x) = \nu\Delta
  \sol(t,x) + \frac{\partial W}{\partial t}(t,x) \\
  &\sol(0,x)=\sol_0(x), 
  \end{aligned}\right.
\end{equation}
where $x=(x_1,x_2) \in \TT^2$, the two-dimensional torus $[0,2\pi]
\times [0,2\pi]$, $\nu >0$ is
the viscosity constant, $\frac{\partial W}{\partial t}$ is a
white-in-time stochastic forcing to be specified below, and 
\begin{align*}
  B(\sol,\tilde \sol)(x) = \sum_{i=1}^2 (\mathcal K\sol)_i(x) \frac{\partial \tilde
  \sol}{\partial x_i}(x)\;,
\end{align*}
where $\mathcal{K}$ is the Biot-Savart
integral operator which will be defined next. First, we define a
convenient basis in which we will perform all explicit calculations.
Setting $\ZZ^2_+=\{(j_1,j_2) \in \ZZ^2 : j_2 > 0\} \cup \{(j_1,j_2) \in
\ZZ^2 : j_1 > 0, j_2=0\}$, $\ZZ^2_-=-\ZZ^2_+$ and $\ZZ^2_0=\ZZ^2_+ \cup
\ZZ^2_-$, we define a real Fourier basis for functions on $\TT^2$
with zero spatial mean by
\begin{align*}
  \base_k(x)=
  \begin{cases}
    \sin(k\cdot x) & k \in \ZZ^2_+\\
    \cos(k\cdot x) & k \in \ZZ^2_-\ .
  \end{cases}
\end{align*}
Write $\sol(t,x)=\sum_{k \in \ZZ^2_0} \solf_k(t) \base_k(x)$ for the
expansion of the solution in this basis.  With this notation, in the
two-dimensional periodic setting,
\begin{align}
  \label{eq:BiotSavart}
  \mathcal{K}(\sol)= \sum_{k \in \ZZ^2_0} \frac{k^\perp}{|k|^2}\solf_k \base_{-k}, 
\end{align}
where $k^\perp=(-k_2,k_1)$.  See for example \cite{b:MajdaBertozzi02}
for more details on the deterministic vorticity formulation in a
periodic domain. We use the vorticity formulation for simplicity, but all
of our results can easily be translated into statements about the velocity
formulation of the problem. We  solve \eqref{eq:vorticity} on the
space $\LL^2=\{ f=\sum_{k\in \ZZ_0^2} a_k e_k : \sum |a_k|^2 < \infty\}$.
For $f=\sum_{k\in \ZZ_0^2} a_k e_k $, we define the norms
$\|f\|^2=\sum |a_k|^2 $ and $\|f\|_1^2=\sum |k|^2|a_k|^2$.

The emphasis of this note will be on forcing which directly excites
only a few degrees of freedom. Such forcing is both of primary
modeling interest and is technically the most difficult.
Specifically we consider forcing of the form
\begin{align}
  \label{eq:W}
  W(t,x)=\sum_{ k \in \mathcal{Z}_*} \sigma_k W_k(t) \base_k(x) \ .
\end{align}
Here $\mathcal{Z}_*$ is a finite subset of $\ZZ^2_0$, $\sigma_k >0$, and $\{ W_k: k
\in \mathcal{Z}_*\}$ is a collection of mutually independent standard
scalar Brownian Motions on a probability space $(\Omega,\mathcal{F},\PP)$. 

The note describes two sets of results contained in two papers by two
different subsets of the authors. In the first, Mattingly and Pardoux
\cite{b:MattinglyPardoux03Pre} give conditions ensuring that any
projection of the time $t$ transition probability of the solution of
\eqref{eq:vorticity} onto a finite dimensional subspace has a
$C^\infty$ density with respect to Lebesgue measure. The result is
based on the Malliavin calculus.  Under additional conditions, this
density is shown to be everywhere positive. The techniques developed
are quite general and we expect they can be applied to many nonlinear,
stochastic partial differential equations with additive noise. These
results provide a first step towards a truly infinite dimensional
version of H\"ormanders celebrated ``sum of squares'' theorem
\cite{b:Hormander85}.

In the second paper, Hairer and Mattingly \cite{b:HairerMattingly04}
give necessary and sufficient conditions for the main results and
estimates of \cite{b:MattinglyPardoux03Pre} to hold. They then use
these tools to build a theory which, when applied to
\eqref{eq:vorticity}, proves that it has a unique invariant measure
under extremely general and essentially sharp assumptions. In addition
to the tools from \cite{b:MattinglyPardoux03Pre}, they introduce new
concept and tool which together provide an abstract framework in which
the ergodicity of \eqref{eq:vorticity} is proven. The concept is a
generalization of the strong Feller property for a Markov process
which, for reasons that will be made clear below, is called the
\textit{asymptotic strong Feller} property. The main feature of this
property is that a diffusion which is irreducible and asymptotically
strong Feller can have at most one invariant measure. It thus yields a
natural generalization of Doob's theorem.  The tool is an approximate
integration by parts formula, in the sense of Malliavin calculus,
which is used to prove that the system enjoys the asymptotic strong
Feller property.  To the best of the authors knowledge, this paper is
the first to prove ergodicity of a nonlinear stochastic partial
differential equation (SPDE) under assumptions comparable to those
assumed when studying finite dimensional stochastic differential
equations.

The ergodic theory of infinite dimensional stochastic systems, and
SPDEs specifically, has been a topic of intense study over the last
two decades. Until recently, the forcing was always assumed to be
elliptic and spatially rough. In our context this translates to $
\mathcal{Z}_*= \ZZ^2_0$ and $|\sigma_k| \sim |k|^{-\alpha}$ for some
positive $\alpha$.  Flandoli and Maslowski \cite{b:FlMa95} first
proved ergodic results for \eqref{eq:vorticity} under such
assumptions. This line of inquiry was extended and simplified in
\cite{b:Fe97,b:GoldysMaslowski04bPre}. They represent a larger body of
literature which characterizes the extent to which classical ideas
developed for finite dimensional Markov processes apply to infinite
dimensional processes. Principally they use tools from infinite
dimensional stochastic analysis to prove that the processes are strong
Feller in an appropriate topology and then deduce ergodicity.

Next three groups of authors in
\cite{b:KuksinShirikyan00,b:BricmontKupiainenLefevere01,b:EMattinglySinai00},
contemporaneously greatly expanded the cases known to be ergodic. They
use the Foias-Prodi type reduction, first adapted to the stochastic
setting in \cite{b:Mattingly98b} and the pathwise contraction of the
high spatial frequencies already used in \cite{b:Mattingly98} to prove
ergodicity of \eqref{eq:vorticity} at sufficiently high viscosity. All
of the results hinged on the observation that if all of the unstable
directions are stochastically perturbed, then the system could be
shown to be ergodic. A general overview of these ideas with simple
examples can be found in \cite{b:Mattingly03Pre}. These ideas have
been continued in a number of papers. See for instance
\cite{b:ELui02,b:BricmontKupiainenLefevere02,b:Mattingly02,b:Hairer02,b:KuksinShirikyan02a,b:Mattingly03Pre}.

Unfortunately, the best current estimates on the number of unstable
directions in \eqref{eq:vorticity} grow inversely with the viscosity
$\nu$. Hence the physically important limit of $\nu \rightarrow 0$
while a fixed, finite scale is forced were previously outside the
scope of the theory. However there existed strong indications that
ergodicity held in this case. Specifically in \cite{b:EMattingly00} it
was shown that the generator of the diffusion associated to finite
dimensional Galerkin approximations of \eqref{eq:vorticity} was
hypoelliptic in the sense of H\"ormander when only a few directions
were forced. This hypoellipticity is the crucial ingredient in the
proof of ergodicity from \cite{b:EMattingly00}.

The ``correct'' ergodic theorem needs to incorporate in its statement
information on how the randomness spreads from the few forced
directions to all of the unstable directions. This understanding when
combined with what had been learned in
\cite{b:Mattingly98b,b:Mattingly98,b:KuksinShirikyan00,b:BricmontKupiainenLefevere01,b:EMattinglySinai00}
should yield unique ergodicity.  This is the program executed in the
papers discussed in this note.

\section{The Geometry of the Forcing and Cascade of Randomness}
\label{sec:conditions}
The geometry of the forcing is encoded in the structure of
$\mathcal{Z}_*$ from \eqref{eq:W}. As observed in
\cite{b:EMattingly00}, its structure gives information about
how the randomness is spread throughout phase space by the
nonlinearity.

Define $\mathcal{Z}_0$ to be the symmetric, and hence translationally
stationary part of the forcing set $\Zf$, given by $\mathcal{Z}_0= \Zf
\cap (-\Zf)$. Then define the collection
\begin{equation*}
  \mathcal{Z}_n= \big\{ \ \ell + j \in \ZZ^2_0 \ :\ 
  j \in \mathcal{Z}_0,\ \ell\in
  \mathcal{Z}_{n-1}
 \quad\text{with}\quad \ell^\perp\cdot j \not =0,\ |j|\not =|\ell|\ \big\}
\end{equation*}
and lastly, 
\begin{equation*}
  \mathcal{Z}_\infty=\bigcup_{n=1}^\infty \mathcal{Z}_n .
\end{equation*}

$\mathcal{Z}_\infty$ captures the directions to which the randomness
has spread. This can be understood in the following way. Denote by
$\d_k$ the partial derivative into the direction $e_k$ of the phase
space and define (on a formal level) the first order differential
operator $\mathcal{X}$ by
\begin{equation*}
\mathcal{X} = \sum_{k\in \ZZ_0^2} \bigl( B(w,w)_k - \nu |k|^2 \bigr) \d_k\;.
\end{equation*}
Then the generator of the Markov process associated to
\eqref{eq:vorticity} is formally given by
\begin{equation*}
\mathcal L = \mathcal X + {1\over 2} \sum_{k \in \mathcal{Z}_*} \sigma_k \d_k^2\;.
\end{equation*}
Note that $B(w,w)_k = \sum_{\ell,j} c_{k,j,\ell} w_\ell w_j$, where
$c_{k,j,\ell} \neq 0$ if and only if $k \in \{j\pm \ell, -j\pm\ell\}$
and $\ell^\perp\cdot j \not =0$, $|j|\not =|\ell|$.  Therefore, all
differential operators of the type $\d_k$ with $k \in \mathcal
Z_\infty$ can be obtained as an iterated Lie bracket of finite length
involving $X_0$ and $\d_\ell$ with $\ell \in \mathcal Z_*$.

Since we want to ensure that all of the unstable
directions are stochastically agitated, we seek conditions where
$\mathcal{Z}_\infty=\ZZ^2_0$. The following essentially sharp
characterization of this situation is given in \cite{b:HairerMattingly04}.
\begin{proposition}\label{propZ} One has  $\mathcal{Z}_\infty=\ZZ^2_0$ if and only
  if both:
  \begin{enumerate}
  \item Integer linear combinations of elements of $\mathcal{Z}_0$
  generate $\ZZ^2_0$.
\item There exist at least two elements in  $\mathcal{Z}_0$  with
  unequal euclidean norm.
  \end{enumerate}
\end{proposition}

This characterization is sharp in the sense that if $\mathcal Z_* =
-\mathcal Z_*$ and one of the above two conditions fails, then there
exists a non-trivial subspace of $\LL^2$ which is left invariant under
the dynamics of \eqref{eq:vorticity}. Also notice that if
\begin{equation*}
\mathcal{Z}_0= \{ (0,1),(0,-1),(1,1),(-1,-1)\}  
\end{equation*}
then Proposition \ref{propZ} implies that
$\mathcal{Z}_\infty=\ZZ^2_0$. Hence forcing four well chosen modes is
sufficient to have the randomness move through the entire system. Of
course one can also force a small number of modes center elsewhere
than at the origin and obtain the same effect. The next two sections
discuss the implications of $\mathcal Z_\infty =\ZZ^2_0$.

\section{Malliavin Calculus and Densities}
\label{sec:Mal}

We define 
 \begin{equation}
   \label{eq:S}
   S_\infty=\mathrm{Span}\bigl( \base_k: k \in \mathcal{Z}_\infty \cup
   \Zf\bigr)\ . 
\end{equation}
One of the main results of \cite{b:MattinglyPardoux03Pre} is the following:
\begin{theorem}\label{thm:main}
  For any $t > 0$ and any finite dimensional subspace $S$ of $S_\infty$,
  the law of the orthogonal projection $\Pi \sol(t, \cdot)$ of
  $\sol(t, \cdot)$ onto $S$ is absolutely continuous with respect to
  the Lebesgue measure on $S$ and has a $C^\infty$ density.
\end{theorem}

In \cite{b:EckmannHairer00}, Eckmann and Hairer used Malliavin
calculus to prove a version of H\"ormander's ``sum of squares''
theorem for a particular SPDE and deduce ergodicity. However, all of
the techniques of that paper required that the forcing excite all but
a finite number of directions and that the forcing be spatially rough
as in \cite{b:FlMa95,b:DaZa96}. The proof of Theorem \ref{thm:main}
builds on ideas introduced into Malliavin calculus by Ocone in
\cite{b:Ocone88}. The central idea is an alternative representation of
the Malliavin matrix of \eqref{eq:vorticity} using the time reversed
adjoint of the linearization of \eqref{eq:vorticity}. Ocone used this
representation when the SPDE was linear in the initial data and the
forcing. When the noise is additive, \cite{b:MattinglyPardoux03Pre}
extends that idea to the nonlinear case.

Let $J_{s,t}\xi$ be the solution of linearization of \eqref{eq:vorticity}
at time $t$ with initial condition $\xi$ at time $s$, $s \leq t$. Let $\bar
J^*_{s,t}\xi$ denote the solution to the $\LL^2$-adjoint of the
linearizion at time $s$, $s \leq t$, with terminal condition $\xi$ at time
$t$. Since the equation is time reversed, the adjoint is well
posed. With this notation, the so--called  ``Malliavin
covariance matrix'' $\mathcal{M}_t$ can be represented by
\begin{align*}
\langle \mathcal{M}_t \phi ,\phi\rangle = \sum_{k \in
  \mathcal{Z}_*}\int_0^t \sigma_k^2\langle  J_{s,t}e_k, \phi\rangle^2\,ds
  =  \sum_{k \in
  \mathcal{Z}_*}\int_0^t \sigma_k^2\langle  e_k, \bar
  J_{s,t}^*\phi\rangle^2\,ds
\end{align*}
where $\phi \in \LL^2$.
The second of these representations is the one used in
\cite{b:MattinglyPardoux03Pre}.  Because of the time reversal, the
representation is not adapted to the filtration generated by $W$ and
new estimates concerning anticipating stochastic processes are
required to obtain the needed estimates. Essentially one needs to show
that the Malliavin matrix is non-degenerate on the subspace $S$ and
that the moments of the reciprocal of the norm of the Malliavin matrix on this
subspace are finite. This is accomplished through the
following estimate which also gives information about the separation of
the randomness on large and small scales.
\begin{theorem} \label{l:onePhi} Let $\Pi$ be the orthogonal projection of
  $\LL^2$ onto a
  finite dimensional subspace of $S_\infty$. For any $t>0$, $\eta >0$,
  $p \geq 1$, $M >0$ and $K \in (0,1)$ there exist two constants
  $c=c(\nu,\eta,p,|\Zf|,t,K,M,\Pi)$ and 
  $\epsilon_0=\epsilon_0(\nu,K,|\Zf|,t,M,\Pi)$ such that for all $\epsilon \in
  (0,\epsilon_0]$,
  \begin{align*}
    \PP\Bigl( \inf_{\phi \in S(M,K,\Pi)}
    \langle\mathcal{M}_t\phi,\phi\rangle < \epsilon \Bigr) \leq
    c\exp(\eta \|\sol(0)\|^2) \epsilon^p
  \end{align*}
where $S(M,K,\Pi)=\{ \phi \in S_\infty: \|\phi\|=1, \| \phi \|_1\leq M, \| \Pi
\phi\| \geq K \}$.
\end{theorem}

With additional assumptions on the controllability of
\eqref{eq:vorticity} conditions are also given ensuring the strict
positivity of the density. This extends results of Ben Arous and
L\'eandre \cite{b:BenArousLeandre91b} and Aida, Kusuoka and Stroock
\cite{b:AidaKusuokaStroock93} to this setting. We refer the reader to
\cite{b:MattinglyPardoux03Pre} for the exact conditions and the
details.

\section{Unique Ergodicity}
\label{sec:erg}
Recall that an \textit{invariant measure} for \eqref{eq:vorticity} is
a probability measure $\mu_\star$ on $\LL^2$ such that $P_t^*\mu_\star
= \mu_\star$, where $P_t^*$ is the semigroup on measures dual to the
Markov transition semigroup $P_t$ defined by
$(P_t\phi)(\sol)=\EE_{\sol}\phi(\sol_t)$ with $\phi\in C_b(\LL^2)$.
While the existence of an invariant measure for \eqref{eq:vorticity}
can be proved by ``soft'' techniques using the regularizing and
dissipativity properties of the flow \cite{b:Cruzeiro90,b:Fl94},
showing its uniqueness is a more challenging problem that requires a
detailed analysis of the nonlinearity. The importance of showing the
uniqueness of $\mu_\star$ is illustrated by the fact that it implies
that
\begin{equation}
\lim_{T \to \infty} \frac{1}{T} \EE \int_0^T \phi(w_t)\,dt = \int_{\LL^2} \phi(w)\,\mu_\star(dw)\;,
\end{equation}
for all bounded continuous functions $\phi$ and all initial conditions
$w_0 \in \LL^2$. It thus gives some mathematical ground to the
\textit{ergodic assumption} usually made in the physics literature
when discussing the qualitative behavior of \eqref{eq:vorticity}.  The
following theorem is the main result of \cite{b:HairerMattingly04}.

\begin{theorem}\label{theo:realmain}
  If $\mathcal{Z}_\infty=\ZZ^2_0$, then, \eqref{eq:vorticity} has a
  unique invariant measure in $\LL^2$.
\end{theorem}
When combined with Proposition \ref{propZ}, this theorem gives easy to
verify conditions guaranteeing a unique invariant measure.

The concept of a strong Feller Markov process appears to be less
useful in infinite dimensions than in finite dimensions. In
particular if $P_t$ is strong Feller, then the measures
$P_t(u,\cdot)$ and $P_t(v,\cdot)$ are equivalent for all initial
conditions $u,v \in \LL^2$. It is easy to construct an ergodic SPDE which
does not satisfy this property.

Recall the following standard sufficient criteria for $P_t$ to be strong Feller~:
 there exists a locally bounded function $C(\sol,t)$ such that
\begin{align*}
|\nabla  (P_t \phi)(\sol)| \leq  C(\sol,t)\|\phi\|_\infty
\end{align*}
for all Fr\'echet differentiable functions $\phi:\LL^2 \rightarrow
\RR$. While we will not give the exact definition of the asymptotic
strong Feller property here, the following similar condition implies
that the process is asymptotically strong Feller: there exists a
locally bounded $C(\sol)$, a non-decreasing sequence of times $t_n$,
and a strictly decreasing sequence $\epsilon_n$ with $\epsilon_n
\rightarrow 0$ so that
\begin{align}\label{asf}
|\nabla  (P_{t_n} \phi)(\sol)| \leq  C(\sol)\|\phi\|_\infty + \epsilon_n \|\nabla\phi\|_\infty
\end{align}
for all Fr\'echet differentiable functions $\phi:\LL^2\rightarrow \RR$ and all $n\geq 1$. In applications one
typically has $t_n \rightarrow \infty$. Hence,  the process behaves as if it acquired
the strong Feller property at time infinity, which justifies the
term asymptotic strong Feller.

First observe that $\langle\nabla_{w} (P_{t}
\phi)(w),\xi\rangle=\EE_{w}(\nabla\phi)(\sol_t)J_{0,t}\xi$. Next we seek a
direction $v$ in the Cameron-Martin space so that if $\mathcal{D}^v$
denotes the Malliavin derivative in the direction $v$ then
$J_{0,t}\xi= \mathcal{D}^v\sol_t$. In finite dimensions, we can often
do this exactly; however, in infinite dimensions we only  know how to achieve this up to
some error. Setting $\rho_t=J_{0,t}\xi - \mathcal{D}^v\sol_t$, we have
the approximate integration by parts formula.
\begin{align*}
  \EE_{w}(\nabla\phi)(\sol_t)J_{0,t}\xi &=
  \EE_{w}\mathcal{D}^v[\phi(\sol_t)] +
  \EE_{w}(\nabla\phi)(\sol_t)\rho_t\\
  &=\EE_{w}\phi(\sol_t)\int_0^t v_s dW_s +
  \EE_{w}(\nabla\phi)(\sol_t)\rho_t\ .
\end{align*}
From this equality one can quickly deduce \eqref{asf}, provided $\EE|\int_0^\infty v_s
dW_s| < \infty$ and $\EE|\rho_t| \rightarrow 0$ as $t \rightarrow
\infty$. In \cite{b:HairerMattingly04}, a $v_t$ is chosen so that
these conditions hold. The analysis is complicated by the fact that
the $v_t$ constructed there is not adapted to the Brownian filtration. This complication
seems unavoidable. Hence, the stochastic integral is a Skorohod integral
and all of the calculations are made more complicated.
 
The ideas developed here can also be used to prove exponential
mixing using the ideas from \cite{b:Mattingly02,b:Hairer02}. These
results will be presented elsewhere.


\begin{thebibliography}{10}

\bibitem{b:AidaKusuokaStroock93}
S.~Aida, S.~Kusuoka, and D.~Stroock.
\newblock On the support of {W}iener functionals.
\newblock In {\em Asymptotic problems in probability theory: Wiener functionals
  and asymptotics (Sanda/Kyoto, 1990)}, volume 284 of {\em Pitman Res. Notes
  Math. Ser.}, pages 3--34. Longman Sci. Tech., Harlow, 1993.

\bibitem{b:BenArousLeandre91b}
G.~Ben~Arous and R.~L{\'e}andre.
\newblock D\'ecroissance exponentielle du noyau de la chaleur sur la diagonale.
  {II}.
\newblock {\em Probab. Theory Related Fields}, 90(3):377--402, 1991.

\bibitem{b:BricmontKupiainenLefevere01}
J.~Bricmont, A.~Kupiainen, and R.~Lefevere.
\newblock Ergodicity of the 2{D} {N}avier-{S}tokes equations with random
  forcing.
\newblock {\em Comm. Math. Phys.}, 224(1):65--81, 2001.
\newblock Dedicated to Joel L. Lebowitz.

\bibitem{b:BricmontKupiainenLefevere02}
J.~Bricmont, A.~Kupiainen, and R.~Lefevere.
\newblock Exponential mixing of the 2{D} stochastic {N}avier-{S}tokes dynamics.
\newblock {\em Comm. Math. Phys.}, 230(1):87--132, 2002.

\bibitem{b:Cruzeiro90}
Ana~Bela Cruzeiro.
\newblock Solutions et mesures invariantes pour des \'equations d'\'evolution
  stochastiques du type {N}avier-{S}tokes.
\newblock {\em Exposition. Math.}, 7(1):73--82, 1989.

\bibitem{b:DaZa96}
Giuseppe Da~Prato and Jerzy Zabczyk.
\newblock {\em Ergodicity for Infinite Dimensional Systems}.
\newblock Cambridge, 1996.

\bibitem{b:ELui02}
Weinan E and Di~Liu.
\newblock Gibbsian dynamics and invariant measures for stochastic dissipative
  {PDE}s.
\newblock {\em Journal of Statistical Physics}, 108(5/6):1125--1156, 2002.

\bibitem{b:EMattinglySinai00}
Weinan E, J.~C. Mattingly, and Ya~G. Sinai.
\newblock Gibbsian dynamics and ergodicity for the stochastic forced
  {N}avier-{S}tokes equation.
\newblock {\em Comm. Math. Phys.}, 224(1), 2001.

\bibitem{b:EMattingly00}
Weinan E and Jonathan~C. Mattingly.
\newblock Ergodicity for the {N}avier-{S}tokes equation with degenerate random
  forcing: finite-dimensional approximation.
\newblock {\em Comm. Pure Appl. Math.}, 54(11):1386--1402, 2001.

\bibitem{b:EckmannHairer00}
J.-P. Eckmann and M.~Hairer.
\newblock Uniqueness of the invariant measure for a stochastic {PDE} driven by
  degenerate noise.
\newblock {\em Comm. Math. Phys.}, 219(3):523--565, 2001.

\bibitem{b:Fe97}
Benedetta Ferrario.
\newblock Ergodic results for stochastic {N}avier-{S}tokes equation.
\newblock {\em Stochastics and Stochastics Reports}, 60(3--4):271--288, 1997.

\bibitem{b:Fl94}
Franco Flandoli.
\newblock Dissipativity and invariant measures for stochastic {N}avier-{S}tokes
  equations.
\newblock {\em NoDEA}, 1:403--426, 1994.

\bibitem{b:FlMa95}
Franco Flandoli and B.~Maslowski.
\newblock Ergodicity of the {2-D} {N}avier-{S}tokes equation under random
  perturbations.
\newblock {\em Comm. in Math. Phys.}, 171:119--141, 1995.

\bibitem{b:GoldysMaslowski04bPre}
B.~Goldys and B.~Maslowski.
\newblock Exponential ergodicity for stochastic burgers and {2D}
  {N}avier-{S}tokes equations.
\newblock Preprint, August 2004.

\bibitem{b:Hairer02}
Martin Hairer.
\newblock Exponential mixing properties of stochastic {PDE}s through asymptotic
  coupling.
\newblock {\em Probab. Theory Related Fields}, 124(3):345--380, 2002.

\bibitem{b:HairerMattingly04}
Martin Hairer and Jonathan~C. Mattingly.
\newblock Ergodicity of the degenerate stochstic {2D} {N}avier--{S}tokes
  equation.
\newblock Submitted, June 2004.

\bibitem{b:Hormander85}
Lars H{\"o}rmander.
\newblock {\em The Analysis of Linear Partial Differential Operators {I--IV},
  publisher = {Springer}, address = {New York}, year = 1985}.

\bibitem{b:KuksinShirikyan00}
Sergei Kuksin and Armen Shirikyan.
\newblock Stochastic dissipative {P}{D}{E}s and {G}ibbs measures.
\newblock {\em Comm. Math. Phys.}, 213(2):291--330, 2000.

\bibitem{b:KuksinShirikyan02a}
Sergei Kuksin and Armen Shirikyan.
\newblock Coupling approach to white-forced nonlinear {PDE}s.
\newblock {\em J. Math. Pures Appl. (9)}, 81(6):567--602, 2002.

\bibitem{b:MajdaBertozzi02}
Andrew~J. Majda and Andrea~L. Bertozzi.
\newblock {\em Vorticity and incompressible flow}, volume~27 of {\em Cambridge
  Texts in Applied Mathematics}.
\newblock Cambridge University Press, Cambridge, 2002.

\bibitem{b:Mattingly98b}
Jonathan~C. Mattingly.
\newblock {\em The Stochastically forced {N}avier-{S}tokes equations: energy
  estimates and phase space contraction}.
\newblock PhD thesis, Princeton University, 1998.

\bibitem{b:Mattingly98}
Jonathan~C. Mattingly.
\newblock Ergodicity of $2${D} {N}avier-{S}tokes equations with random forcing
  and large viscosity.
\newblock {\em Comm. Math. Phys.}, 206(2):273--288, 1999.

\bibitem{b:Mattingly02}
Jonathan~C. Mattingly.
\newblock Exponential convergence for the stochastically forced
  {N}avier-{S}tokes equations and other partially dissipative dynamics.
\newblock {\em Comm. Math. Phys.}, 230(3):421--462, 2002.

\bibitem{b:Mattingly03Pre}
Jonathan~C. Mattingly.
\newblock On recent progress for the stochastic {N}avier {S}tokes equations.
\newblock In {\em Journ\'ees \'Equations aux d\'eriv\'ees partielles},
  Forges-les-Eaux, 2003.
\newblock see http://www.math.sciences.univ-nantes.fr/edpa/2003/html/.

\bibitem{b:MattinglyPardoux03Pre}
Jonathan~C. Mattingly and \'Etienne Pardoux.
\newblock Malliavin calculus and the randomly forced {N}avier {S}tokes
  equation.
\newblock Submitted to Jounal of the MAS, 2004.

\bibitem{b:Ocone88}
Daniel Ocone.
\newblock Stochastic calculus of variations for stochastic partial differential
  equations.
\newblock {\em J. Funct. Anal.}, 79(2):288--331, 1988.

\end{thebibliography}
\end{document}